\documentclass[11pt]{article}
\usepackage{geometry}                
\usepackage[parfill]{parskip}    
\usepackage{graphicx}
\usepackage{amssymb}
\usepackage{epstopdf}
\usepackage{amsmath}
\usepackage{mathtools}
\usepackage[english]{babel}
\usepackage{amsthm}
\usepackage[usenames,dvipsnames]{xcolor}
\usepackage{amscd}
\usepackage{scalerel}
\usepackage{soul}

\newtheorem{theorem}{Theorem}

\newtheorem{corollary}[theorem]{Corollary}

\newcommand{\DQ}{\operatorname{DQ}}
\newcommand{\ds}{\displaystyle}

\title{Recognizing difference quotients of real functions.}
\author{Trevor Richards, Jimmy Yau\footnote{This paper represents the content of the second author's undergraduate research project while at Washington and Lee University in the Fall term of 2015 and the Winter term of 2016, under the direction of the first author.}}

\begin{document}
\maketitle

\abstract{
\noindent For a real function $f:[0,1]\to\mathbb{R}$, the difference quotient of $f$ is the function of two real variables $\DQ_f(a,b)=\dfrac{f(b)-f(a)}{b-a}$, which we view as defined on the triangle $\mathcal{T}=\{(a,b):0\leq a<b\leq1\}$.  In this paper we investigate how to determine whether a given function of two variables $H(a,b)$ is the difference quotient of some real function $f(x)$.  We develop three independent methods for recognizing such a function $H$ as a difference quotient, and corresponding methods for recovering the underlying function $f$ from $H$.}

\section{Introduction}
Given a function $f(x)$ defined at two distinct points $a,b\in[0,1]$, the difference quotient of $f$ from $a$ to $b$ is the slope of the line segment connecting points $(a,f(a))$ and $(b,f(b))$. This we denote as \[\DQ_f(a,b)=\dfrac{f(b)-f(a)}{b-a}.\]

However, simplification and rearrangement of $\DQ_f$ may leave this function in a form that is unrecognizable as a difference quotient. Therefore, it is desirable to develop some means to determine whether a given function of two variables $H(a,b)$ is the difference quotient of some function $f(x)$ and, if so, a method for recovering $f$ from $H$.  We will give three such tests for whether $H$ is equal to $\DQ_f$ for some $f$, and also corresponding methods to recover $f$ from $H$.

Throughout this paper, all variables are taken to be real, and all functions are assumed to take real values.  We will restrict our attention to functions $f(x)$ having domain $[0,1]$.  Therefore we will consider $H$ to have as its domain either the square $[0,1]\times[0,1]$, or some subset thereof.  Of course, if the difference quotient $\DQ_f$ extends continuously to the diagonal \[\mathcal{D}=\{(a,a):a\in[0,1]\},\] then $f$ must be differentiable and $\DQ_f=f'$ on $\mathcal{D}$.

The tests we provide are presented in an increasing order of specificity with regard to the sort of functions $H$ to which they may be applied, and along the way, we will give concrete examples of functions $H$ which we may test.

The first criteria, which we call the \textit{Algebraic Criteria} and prove in Section~\ref{sect: Proof of Algebraic Criteria.}, is meant to detect the difference quotient of a completely arbitrary real function $f(x)$. Therefore the function $H$ is not expected to be defined on $\mathcal{D}$.  It must, however, be symmetric (ie $H(a,b)=H(b,a)$) in order to be a difference quotient, so, rather than include symmetry as a part of the criteria, we merely assume that $H$ is defined only on the upper triangle \[\mathcal{T}=\{(a,b):0\leq a<b\leq1\}.\]  If we wish to extend $H$ to the lower triangle we may do so symmetrically by defining $H(b,a)=H(a,b)$ for $a<b$.

\noindent\rule{\textwidth}{0.1pt}
\begin{theorem}[Algebraic Criteria]\label{Algebraic Criteria}

Let $H: \mathcal{T}\rightarrow \mathbb{R}$ be given. Then the following are equivalent:
\bigskip
\begin{enumerate}
\item There exists some $f:[0,1]\rightarrow\mathbb{R}$ such that $\DQ_f(a,b)=H(a,b)$.
\item For all $0\leq a<b<c\leq1$,
\[H(a,c)=\frac{H(a,b)(b-a)+H(b,c)(c-b)}{c-a}.\]

\item For all $0<b<c\leq1$,
\[H(0,c)=\frac{H(0,b)(b)+H(b,c)(c-b)}{c}.\]

\end{enumerate}
\end{theorem}
\noindent\rule{\textwidth}{0.1pt}

In the proof of Theorem~\ref{Algebraic Criteria}, which we give in Section~\ref{sect: Proof of Algebraic Criteria.}, we construct the underlying function $f$ assuming $H$ satisfies the second and third items of Theorem~\ref{Algebraic Criteria}. We enshrine that construction in the following corollary.

\begin{corollary}\label{Algebraic Corollary}
For $H(a,b)$ satisfying Items 2 and 3 of Theorem~\ref{Algebraic Criteria}, the functions $f(x)$ which make $H=\DQ_f$ are all functions of the form \[f(x)=xH(0,x)+C\] for a constant $C$.
\end{corollary}

We observe furthermore that the expression
\[H(a,c)=\frac{H(a,b)(b-a)+H(b,c)(c-b)}{c-a}\] from Item~2 of Theorem~\ref{Algebraic Criteria} may be rewritten as
\[H(b,c)(b-c)-H(a,c)(a-c)+H(a,b)(a-b)=0.\]  The left hand side of this equation may easily be recognized as the determinant of the matrix \[M=\begin{bmatrix}
       H(b,c) & H(a,c) & H(a,b)\\[0.3em]
       a & b & c \\[0.3em]
       1 & 1 & 1
     \end{bmatrix},\]
which will be referred to in Corollary~\ref{cor: Algebraic Criteria Matrix} below.  From this observation, we extract several secondary tests, which follow from Theorem~\ref{Algebraic Criteria}.  It is interesting to point out that, in the case $H=\DQ_f$ for some function $f$, the dimension and null set of $M$ are in fact independent of $H$.

\begin{corollary}\label{cor: Algebraic Criteria Matrix}
Let $H: \mathcal{T}\rightarrow \mathbb{R}$ be given.  The following are equivalent:

\bigskip
\begin{enumerate}
\item There exists some $f:[0,1]\rightarrow\mathbb{R}$ such that $\DQ_f(a,b)=H(a,b)$.
\item For all $0\leq a<b<c\leq0$, $\det(M)=0$.
\item For all $0\leq a<b<c\leq0$, $\dim(M)=2$.
\item For all $0\leq a<b<c\leq0$, $\operatorname{Null}(M)=\left\{k\begin{bmatrix}
       \frac{b-c}{c-a}\\[0.3em]
       1\\[0.3em]
       \frac{a-b}{c-a}
       \end{bmatrix}
:k\in\mathbb{R}\right\}$.

\end{enumerate}
\end{corollary}

We now turn our attention to the difference quotients of differentiable functions $f(x)$.  If $f$ is differentiable, then its difference quotient $\DQ_f$ may now defined on $\overline{\mathcal{T}}={\mathcal{T}}\cup{\mathcal{D}}$ by \[\DQ_f(a,a)=f'(a).\]  With this definition in mind, we obtain the following test, which we call the \textit{Integrable Criteria} and prove in Section~\ref{sect: Integrable Criteria.}.

\noindent\rule{\textwidth}{0.1pt}

\begin{theorem}[Integrable Criteria]\label{Integrable Criteria}

\noindent Let $H: \overline{\mathcal{T}}\rightarrow \mathbb{R}$ be given, such that the function $h(x)=H(x,x)$ is (Lebesgue) integrable on $[0,1]$. Then the following are equivalent:
\bigskip

\begin{enumerate}
\item There exists some $f:[0,1]\rightarrow\mathbb{R}$ such that $\DQ_f(a,b)=H(a,b)$.

\item For all $0\leq a<b\leq1$, \[(b-a)H(a,b)=\int_a^b H(s,s)ds.\]
\end{enumerate}
\end{theorem}
\noindent\rule{\textwidth}{0.1pt}

As with the Algebraic Criteria, the proof of the Integrable Criteria consists of an explicit construction of the functions $f$ for which $H=\DQ_f$.  This construction is recorded as the following corollary.

\begin{corollary}\label{Integrable Corollary}
For $H(a,b)$ satisfying Item 2 of Theorem~\ref{Integrable Criteria}, the functions $f(x)$ which make $H=\DQ_f$ are all functions of the form \[f(x)=\int_0^x H(s,s)ds+C\] for a constant $C$.
\end{corollary}

\noindent In Section~\ref{sect: Integrable Criteria.}, we will also observe that, for a differentiable function $f$, the sum of the first partial derivatives of the difference quotient of $f$ is the difference quotient of $f'$.

\begin{theorem}
Let $f$ be differentiable on $[0,1]$.  Then \[\left(\dfrac{\partial}{\partial a}+\dfrac{\partial}{\partial b}\right)\DQ_f(a,b)=\DQ_{f'}(a,b).\]
\end{theorem}

Our final test, which we call the \textit{Summation Criteria} and will prove in Section~\ref{sect: Summation Criteria.}, is designed to detect the difference quotient of an analytic function $f$.  The key observation is that the difference quotient of the function $f(x)=x^n$ is exactly \[\DQ_f(a,b)=\ds\sum_{i=0}^{n-1}a^ib^{n-1-i}.\]

\noindent\rule{\textwidth}{0.1pt}
\begin{theorem}[Summation Criteria]\label{Summation Criteria}
Let $H:\overline{\mathcal{T}}\rightarrow \mathbb{R}$ be an analytic function $H(a,b)=\ds\sum_{i,j\geq0}c_{ij}a^ib^j$ on $\overline{\mathcal{T}}$. The following are equivalent:

\bigskip
\begin{enumerate}
\item There exists some $f:[0,1]\rightarrow\mathbb{R}$ such that $\DQ_f(a,b)=H(a,b)$.

\item For each $p\in\mathbb{N},$ there exists a constant $c_p$ such that for all $i,j\geq0$ with $i+j=p$, $c_{ij}=c_p$.
\end{enumerate}
\end{theorem}
\noindent\rule{\textwidth}{0.1pt}

As with the earlier criteria, for functions $H$ that satisfy the items of Theorem~\ref{Summation Criteria}, we obtain an explicit construction of the functions $f$ that will make $H=\DQ_f$.

\begin{corollary}\label{Summation Corollary}

For $H(a,b)$ satisfying Item 2 of Theorem~\ref{Summation Criteria}, the functions $f(x)$ which make $H=\DQ_f$ are all functions of the form \[f(x)=\left(\sum\limits_{p=0}^{\infty}c_px^{p+1}\right)+C\] for a constant $C$.

\end{corollary}

\section{Algebraic Criteria}\label{sect: Proof of Algebraic Criteria.}

The Algebraic Criteria is the most general of the tests we will develop in this paper.  It is designed to detect the difference quotient of an arbitrary function $f:[0,1]\to\mathbb{R}$, and thus, the function $H$ being tested need only be defined on the upper triangle $\mathcal{T}$.

For such a function $H$, we will show that the following items are equivalent.

\medskip

\noindent\textbf{Algebraic Criteria}

\begin{enumerate}
\item\label{item: Algebraic Criteria Item DQ_f.} There exists some $f:[0,1]\rightarrow\mathbb{R}$ such that $\DQ_f(a,b)=H(a,b)$.

\item\label{item: Algebraic Criteria Item all a,b,c.} For all $0\leq a<b<c\leq1$,
\[H(a,c)=\frac{H(a,b)(b-a)+H(b,c)(c-b)}{c-a}.\]

\item\label{item: Algebraic Criteria Item all b,c.} For all $0<b<c\leq1$,
\[H(0,c)=\frac{H(0,b)(b)+H(b,c)(c-b)}{c}.\]
\end{enumerate}


The fact that Item~\ref{item: Algebraic Criteria Item DQ_f.} implies Item~\ref{item: Algebraic Criteria Item all a,b,c.} follows directly from the definition of the difference quotient, and of course Item~\ref{item: Algebraic Criteria Item all a,b,c.} implies Item~\ref{item: Algebraic Criteria Item all b,c.} after setting $a=0$.  Therefore, we proceed to a more interesting fact that Item~\ref{item: Algebraic Criteria Item all b,c.} implies Item~\ref{item: Algebraic Criteria Item DQ_f.}.

For a fixed constant $C\in\mathbb{R}$, we define $f:[0,1]\to\mathbb{R}$ by

$$f(x)=\begin{cases}
C&x=0\\
xH(0,x)+C&x\neq0
\end{cases}.$$

With this choice of $f$, we will show that Item~\ref{item: Algebraic Criteria Item all b,c.} implies $H=\DQ_f$. Note that it is easy to show that two real valued functions have the same difference quotient if and only if they differ by a constant. Therefore showing that $H=\DQ_f$ will immediately establish Corollary~\ref{Algebraic Corollary} as well.

\begin{proof}

Assume that Item~3 of the Algebraic Criteria holds, and let $0\leq i<j\leq1$ be given.  We wish to show that, for the function $f$ defined above, $\DQ_f(i,j)=H(i,j)$.  We have two values of $i$ to check: 

\noindent Case $i=0$: $$\DQ_f(0,j)=\dfrac{f(j)-f(0)}{j-0}=\dfrac{jH(0,j)+C-C}{j}=H(0,j).$$

\noindent Case $0<i$:  Note that the condition present in Item~3 of the Algebraic Criteria may be rearranged to give that, for all $0<b<c\leq1$,
\[H(b,c)=\frac{cH(0,c)-bH(0,b)}{c-b}.\]  We therefore have that $$\DQ_f(i,j)=\dfrac{f(j)-f(i)}{j-i}=\dfrac{(jH(0,j)+C)-(iH(0,i)+C)}{j-i}=\dfrac{jH(0,j)-iH(0,i)}{j-i}=H(i,j).$$  This establishes the desired result.

\end{proof}

\subsection{Application to Linear Algebra}

Continuing with our analysis of an arbitrary function $H$ defined on $\mathcal{T}$, for any $0\leq a<b<c\leq1$, we define the matrix

\[M=\begin{bmatrix}
       H(b,c) & H(a,c) & H(a,b)\\[0.3em]
       a & b & c \\[0.3em]
       1 & 1 & 1
     \end{bmatrix}.\]

We will now show the equivalence of the following items from Corollary~\ref{cor: Algebraic Criteria Matrix}.

\begin{enumerate}
\item\label{Item: Algebraic Criteria Matrix Item DQ_f.} There exists some $f:[0,1]\rightarrow\mathbb{R}$ such that $\DQ_f(a,b)=H(a,b)$.

\item\label{Item: Algebraic Criteria Matrix Item zero determinant.} For all $0\leq a<b<c\leq0$, $\det(M)=0$.

\item\label{Item: Algebraic Criteria Matrix Item dimension two.} For all $0\leq a<b<c\leq0$, $\dim(M)=2$.

\item\label{Item: Algebraic Criteria Matrix Item fixed solution set.} For all $0\leq a<b<c\leq0$, $\operatorname{Null}(M)=\left\{k\begin{bmatrix}
       \frac{b-c}{c-a}\\[0.3em]
       1\\[0.3em]
       \frac{a-b}{c-a}
       \end{bmatrix}
:k\in\mathbb{R}\right\}$.

\end{enumerate}

It may easily be shown that the condition in Item~\ref{item: Algebraic Criteria Item all a,b,c.} of the Algebraic Criteria may be rewritten as the condition $\det(M)=0$.  Therefore, the Algebraic Criteria gives us the equivalence of Items~\ref{Item: Algebraic Criteria Matrix Item DQ_f.}~and~\ref{Item: Algebraic Criteria Matrix Item zero determinant.} from Corollary~\ref{cor: Algebraic Criteria Matrix}.  Moreover, basic linear algebra shows that, among the items of Corollary~\ref{cor: Algebraic Criteria Matrix}, Item~\ref{Item: Algebraic Criteria Matrix Item fixed solution set.} implies Item~\ref{Item: Algebraic Criteria Matrix Item dimension two.}, which in turn implies Item~\ref{Item: Algebraic Criteria Matrix Item zero determinant.}. Thus, we have the following implications.
$$1\Longleftrightarrow2\Longleftarrow3\Longleftarrow4$$

We will establish the remaining implications by showing that Item~\ref{Item: Algebraic Criteria Matrix Item DQ_f.} of Corollary~\ref{cor: Algebraic Criteria Matrix} implies Item~\ref{Item: Algebraic Criteria Matrix Item fixed solution set.} of Corollary~\ref{cor: Algebraic Criteria Matrix}.

\begin{proof}

We begin by assuming that Item~\ref{Item: Algebraic Criteria Matrix Item DQ_f.} holds.  We therefore replace $H$ with $\DQ_f$ in the matrix $M$ and seek the solution set for the following matrix equation.

 \[\begin{bmatrix}
       \frac{f(c)-f(b)}{c-b} & \frac{f(c)-f(a)}{c-a} & \frac{f(b)-f(a)}{b-a}\\[0.3em]
       a & b & c \\[0.3em]
       1 & 1 & 1
     \end{bmatrix} \begin{bmatrix}
       x\\[0.3em]
       y\\[0.3em]
       z
     \end{bmatrix}=\vec{0}
     \]

Pulling this matrix equation apart componentwise, we obtain the following system of scalar equations.

\[\begin{cases} x\left(\frac{f(c)-f(b)}{c-b}\right)+y\left(\frac{f(c)-f(a)}{c-a}\right)+z\left(\frac{f(b)-f(a)}{b-a}\right)=0\\xa+yb+zc=0\\x+y+z=0\end{cases}\]

Applying standard algebraic techniques to the latter two equations, we may solve for each of $x$ and $z$ in terms of $y$, obtaining the following.

$$\begin{array}{l}
z=y\left(\frac{b-a}{a-c}\right)\\
x=y\left(\frac{c-b}{a-c}\right)
\end{array}$$

Thus, $\operatorname{Null}(M)$ is contained in the set $$\left\{k\begin{bmatrix}
       \frac{b-c}{c-a}\\[0.3em]
       1\\[0.3em]
       \frac{a-b}{c-a}
       \end{bmatrix}
:k\in\mathbb{R}\right\}.$$

However, it is easy to check that any vector in this set is also in $\operatorname{Null}(M)$.  This completes our proof.

\end{proof}


\subsection{Example}

For $0\leq a<b\leq1$, define

$$H(a,b)=\begin{cases}
0&a,b\in\mathbb{Q}\\
0&a,b\notin\mathbb{Q}\\
\dfrac{1}{b-a}&a\notin\mathbb{Q},b\in\mathbb{Q}\\
\dfrac{-1}{b-a}&a\in\mathbb{Q},b\notin\mathbb{Q}
\end{cases}.
$$

We will verify that this function $H$ satisfies the third item of the Algebraic Criteria, and is therefore the difference quotient of some function $f$.  We will then use Corollary~\ref{Algebraic Corollary} to find the function $f$.

Let $0<b<c\leq1$ be given.  We have four cases to check.

Case $b,c\in\mathbb{Q}$:
$$\dfrac{H(0,b)(b)+H(b,c)(c-b)}{c}=\dfrac{0\cdot b+0\cdot(c-b)}{c}=0=H(0,c).$$

Case $b,c\notin\mathbb{Q}$:
$$\dfrac{H(0,b)(b)+H(b,c)(c-b)}{c}=\dfrac{\dfrac{1}{b-0}\cdot b+0\cdot(c-b)}{c}=\dfrac{1}{c-0}=H(0,c).$$

Case $b\notin\mathbb{Q},c\in\mathbb{Q}$:
$$\dfrac{H(0,b)(b)+H(b,c)(c-b)}{c}=\dfrac{\dfrac{1}{b-0}\cdot b+\dfrac{-1}{c-b}\cdot(c-b)}{c}=0=H(0,c).$$

Case $b\in\mathbb{Q},c\notin\mathbb{Q}$:
$$\dfrac{H(0,b)(b)+H(b,c)(c-b)}{c}=\dfrac{0\cdot b+\dfrac{-1}{c-b}\cdot(c-b)}{c}=\dfrac{-1}{c-0}=H(0,c).$$

Since $H$ satisfies the Algebraic Criteria, we know that $H$ is the difference quotient of a function $f:[0,1]\to\mathbb{R}$.  We may apply Corollary~\ref{Algebraic Corollary}, choosing $C=1$, to find one such $f$.  That is, $H=\DQ_f$, where

$$f(x)=\begin{cases}
1&x=0\\
xH(0,x)+1&x\neq0
\end{cases}.$$

\noindent Thus, we have $f(0)=1$, and if $x\neq0$ then, by the definition of $H$, we have

$$f(x)=\begin{cases}
x\cdot0+1&x\in\mathbb{Q}\\
x\cdot\dfrac{-1}{x}+1&x\notin\mathbb{Q}\end{cases}.$$

\noindent Combining cases and simplifying, we see that $H$ is the difference quotient of the Dirichlet function

$$
f(x)=\begin{cases}1&x\in\mathbb{Q}\\
0&x\notin\mathbb{Q}
\end{cases}.
$$



\section{Integrable Criteria}\label{sect: Integrable Criteria.}

In this section we extend the definition of the difference quotient of $f$ to the diagonal $\mathcal{D}=\{(a,a):0\leq a\leq1\}$ by $\DQ_f(a,a)=f'(a)$, provided that the derivative exists.  With this definition, the fundamental theorem of calculus now suggests the Integrable Criteria, which may be used to detect whether a function $H$ defined on the closed triangle $$\overline{\mathcal{T}}=\{(a,b):0\leq a\leq b\leq1\}$$ is the difference quotient of a differentiable function.  The Integrable Criteria states that for such a function $H$, the following are equivalent.

\medskip

\noindent\textbf{Integrable Criteria}

\begin{enumerate}
\item\label{item: Integrable Criteria Item DQ_f.} There exists some $f:[0,1]\rightarrow\mathbb{R}$ such that $\DQ_f(a,b)=H(a,b)$ on $\overline{\mathcal{T}}$.

\item\label{item: Integrable Criteria Item all a,b.} The function $h(s)=H(s,s)$ is integrable on $[0,1]$, and for all $0\leq a<b\leq1$, \[(b-a)H(a,b)=\int_a^b H(s,s)ds.\]

\end{enumerate}

\noindent The fact that Item~\ref{item: Integrable Criteria Item DQ_f.} implies Item~\ref{item: Integrable Criteria Item all a,b.} follows directly from the fundamental theorem of calculus.  It remains to prove the reverse implication.

\begin{proof}
Assume that the function $H:\overline{\mathcal{T}}\to\mathbb{R}$ satisfies the assumptions of the second item of the Integrable Criteria.  Define the function $f:[0,1]\to\mathbb{R}$ by $$f(x)=\ds\int_0^xH(s,s)ds.$$  Then, by the additivity of the integral, we have $$\DQ_f(a,b)=\dfrac{\ds\int_0^bH(s,s)ds-\int_0^aH(s,s)ds}{b-a}=\dfrac{1}{b-a}\ds\int_a^bH(s,s)ds.$$

By rearranging the equation in Item~\ref{item: Integrable Criteria Item all a,b.} of the Integrable Criteria, we immediately obtain $\DQ_f(a,b)=H(a,b)$.
\end{proof}

Since all functions having the same difference quotient vary only by a constant, we immediately obtain the conclusion found in Corollary~\ref{Integrable Corollary}.

\subsection{Example}
For $0\leq a < b \leq 1$, define $H(a,b)$ to be the average value of $g(x)=e^{x^2}$ on the interval $[a,b]$.  It therefore makes sense to extend this definition to the diagonal by $H(a,a)=e^{a^2}$, for $0\leq a\leq1$.

We check that Item~\ref{item: Integrable Criteria Item all a,b.} of the Integrable Criteria is satisfied.

Let $0\leq a< b\leq1$ be given.  The basic calculus definition of average value of $g(x)$ on the interval $[a,b]$ is
\[\dfrac{1}{b-a}\int_{a}^{b}g(s)ds.\]

Thus we have

$$(b-a)H(a,b)=\dfrac{b-a}{b-a}\int_{a}^{b}e^{s^{2}}ds=\int_{a}^{b}H(s,s)ds.$$

Finally, according to Corollary~\ref{Integrable Corollary}, we have that one of the functions $f$ which makes $H=\DQ_f$ is

$$f(x)=\ds\int_{0}^xe^{s^2}ds.$$

\subsection{Sum of Partial Derivatives}

In this subsection we make note that the sum of the partial derivatives of $\DQ_f$ is equal to $\DQ_{f'}$ $$\left(\dfrac{\partial}{\partial a}+\dfrac{\partial}{\partial b}\right)\DQ_f(a,b)=\DQ_{f'}(a,b).$$  That is, if we define $D^1([0,1])$ to be the set of all differentiable functions from $[0,1]$ to $\mathbb{R}$, $D^1(\mathcal{T})$ to be the set of functions mapping $\mathcal{T}$ to $\mathbb{R}$ both of whose partial derivatives exist, $F([0,1])$ to be the set of functions from $[0,1]$ to $\mathbb{R}$, and $F(\mathcal{T})$ to be the set of functions from $\mathcal{T}$ to $\mathbb{R}$, then the following diagram commutes.

$$
  \scaleto{%
    \begin{CD}
    D^1([0,1]) @>\frac{d}{dx}>> F([0,1]) \\
    @VDQVV @VVDQV\\
    D^1(\mathcal{T}) @>\frac{\partial}{\partial a}+\frac{\partial}{\partial b}>> F(\mathcal{T})
  \end{CD}
}{100pt}
$$

\medskip

\begin{proof}
Using the quotient rule, we observe that

\[\frac{\partial}{\partial a}\dfrac{f(b)-f(a)}{b-a}=\dfrac{-f'(a)(b-a)-(-1)(f(b)-f(a))}{(b-a)^2}\]

and

\[\frac{\partial}{\partial b}\dfrac{f(b)-f(a)}{b-a}=\frac{f'(b)(b-a)-(f(b)-f(a))}{(b-a)^2}.\]

Thus a bit of arithmetic immediately gives

\[\left(\dfrac{\partial}{\partial a}+\dfrac{\partial}{\partial b}\right)\DQ_f(a,b)=\frac{f'(b)-f'(a)}{b-a}.\]

\end{proof}


\section{Summation Criteria}\label{sect: Summation Criteria.}
Our final test is meant to detect the difference quotient of an analytic function.  It comes from the basic observation that the difference quotient of the function $f(x)=x^n$ is $$\DQ_f(a,b)=\dfrac{b^n-a^n}{b-a}=\ds\sum_{i=0}^{n-1}a^ib^{n-1-i}.$$  We see that this is the sum of all terms of the form $a^ib^j$ such that $i,j\geq0$ and $i+j=n-1$.  This observation gives rise to the Summation Criteria, which states for power series \[H(a,b)=\sum\limits_{i,j \geq0}c_{ij}b^{i}a^{j}\] converging absolutely on the closed triangle $\overline{\mathcal{T}}$, that the following are equivalent.
\medskip

\noindent\textbf{Summation Criteria}

\begin{enumerate}
\item\label{item: Summation Criteria Item DQ_f.} There exists some $f:[0,1]\rightarrow\mathbb{R}$ such that $\DQ_f(a,b)=H(a,b)$ on $\overline{\mathcal{T}}$.

\item\label{item: Summation Criteria Item all c.} For each $p\in\mathbb{N},$ there exists a constant $c_p$ such that for all $i,j\geq0$ with $i+j=p$, $c_{ij}=c_p$.
\end{enumerate}

The fact that Item~\ref{item: Summation Criteria Item DQ_f.} implies Item~\ref{item: Summation Criteria Item all c.} follows immediately from the observation above regarding the difference quotient of a power of $x$.  It therefore remains to show that Item~\ref{item: Summation Criteria Item all c.} implies Item~\ref{item: Summation Criteria Item DQ_f.}.

\begin{proof}
Assume that Item~\ref{item: Summation Criteria Item all c.} of the Summation Criteria holds and define $f:[0,1]\rightarrow\mathbb{R}$ by
\[f(x)=\left(\sum\limits_{p=0}^{\infty}c_px^{p+1}\right).\]

\noindent Then

\[
\DQ_f(a,b)=\dfrac{\sum\limits_{p=0}^{\infty}c_pb^{p+1}-\sum\limits_{p=0}^{\infty}c_pa^{p+1}}{b-a}=\sum\limits_{p=0}^{\infty}c_{p}\dfrac{b^{p+1}-a^{p+1}}{b-a}=\sum\limits_{p=0}^{\infty}\sum\limits_{i=0}^{p}c_{p}a^{i}b^{p-i}.
\]

Observe that for each $p\geq1$, and for each $0\leq i\leq p$, setting $j=p-i$ we have $c_{ij}=c_p$.  Therefore, the latter sum may be rewritten

\[
\DQ_f(a,b)=\sum\limits_{p=0}^{\infty}\left(\sum\limits_{\substack{i,j\geq0 \\ i+j=p}}c_{ij}a^{i}b^{j}\right).
\]

\noindent This is indeed a reordering of the power series

\[
H(a,b)=\ds\sum_{i,j\geq0}^\infty c_{ij}a^{i}b^j.
\]
Since $H$ was assumed to be absolutely convergent on $\mathcal{T}$, reorderings of $H$ are equal to $H$, so we have

$$\DQ_f(a,b)=H(a,b).$$

\end{proof}

\subsection{Example}
For $0 \leq a \leq b \leq 1,$ define
\[H(a,b)=\sum_{i,j\geq0}\frac{1}{(i+j)!}a^{i}b^{j}.\]

Before applying the Summation Criteria to $H$, we must verify that $H$ converges absolutely on the closed triangle $\overline{\mathcal{T}}.$  Since each $c_{ij}$ is positive and $0\leq a\leq b\leq1$, we have

\[\sum_{i,j\geq0}\left|\frac{1}{(i+j)!}a^{i}b^{j}\right|<\sum_{i,j\geq0}\frac{1}{(i+j)!}.\]

\noindent We regroup the summation as follows.
\[\sum_{i,j\geq0}\frac{1}{(i+j)!}=\sum_{p=0}^{\infty}\sum_{\substack{i,j\geq0 \\ i+j=p}}\frac{1}{(i+j)!}.\] Looking at the inner summation, we observe that \[\sum_{\substack{i,j\geq0 \\ i+j=p}}\frac{1}{(i+j)!}=\sum_{i=0}^{p}\frac{1}{p!}=(p+1)\frac{1}{p!}.\] Thus, we rewrite \[\sum_{p=0}^{\infty}\sum_{\substack{i,j\geq0 \\ i+j=p}}\frac{1}{(i+j)!}=\sum_{p=0}^{\infty}\frac{p+1}{p!}.\]
Finally, the ratio test guarantees the convergence of $\ds\sum_{p=0}^{\infty}\frac{p+1}{p!}$, so we conclude that $H$ does converge absolutely on $\overline{\mathcal{T}}$.

We observe that, by setting $c_p=\frac{1}{p!}$ for all $p\geq0$, if $i,j\geq0$ with $i+j=p$, then $c_{ij}=c_p$, so that the second item of the Summation Criteria is satisfied.  Using Corollary~\ref{Summation Corollary}, we find one of the functions $f(x)$ for which $H(a,b)=\DQ_f(a,b)$ is \[f(x)=\sum_{p=0}^{\infty}\dfrac{1}{p!}x^{p+1}=x\sum_{p=0}^{\infty}\dfrac{1}{p!}x^{p}.\]  We recognize the series above as the expansion for $f(x)=xe^x$.

\end{document}